\documentclass{article}
\usepackage{amssymb}

\usepackage{graphicx}
\usepackage{amsmath}

\newtheorem{theorem}{Theorem}

\newtheorem{definition}[theorem]{Definition}

\newtheorem{proposition}[theorem]{Proposition}

\newenvironment{proof}[1][Proof]{\textbf{#1.} }{\ \rule{0.5em}{0.5em}}

\newcommand {\Ber}{\mbox{{\rm Ber}}}

\renewcommand{\O}{{\cal O}}
\newcommand{\M}{{\cal M}}
\newcommand{\F}{{\cal F}}

\begin{document}

\title{On the Berezinian of a moduli space of curves in $\mathbf{P}^{n|n+1}$}
\author{Michael Movshev\thanks{This material is partially based upon work supported by the National Science Foundation under agreement No. DMS-0111298. 
}\\IAS \\Princeton , NJ\\0854, USA\\}
\date{\today}

\maketitle
\begin{abstract}
A supermanifold $\mathbf{P}^{3|4}$ is a target  space for twistor string theory.
In this note we identify a line bundle  of holomorphic volume elements $\Ber \M_g(\mathbf{P}^{3|4})$ defined on the moduli space of curves of genus $g$ in $\mathbf{P}^{3|4}$ with a pullback of a line bundle defined on $\M_g(pt)$. We  also give  some generalizations of this fact.
\end{abstract}

\section{Introduction}
A twistor string theory (TST) was introduced in \cite{Twistor} as a tool for studying perturbative expansion of  amplitudes in $D=4, N=4$ Yang-Mills theory (SYM). It was shown later in \cite{BW} that the full TST contains  conformal supergravity and its connection with SYM is not quite clear. However there is a hope  that a modification of TST  is possible that enables us to relate it to perturbative SYM  theory.

It is desirable to achieve a better understanding of the main ingredients of TST. One of such ingredients is an integration measure for loop diagrams. In \cite{Witten} we read:

{\bf...the proper definition of the integration
measure for loop diagrams in twistor space (that is, for
$D$-instanton configurations of positive genus) remains unclear.}

In this note we clarify a situation with the measure.

In TST the target manifold is a complex projective superspace $\mathbf{P}^{3|4}$. We will work in a slightly more general context, when the target is $\mathbf{P}^{n|n+1}$.

Suppose $V$ is equal to $\mathbb{C}^{n+1}$. We denote by  $\Pi$  an operation of parity reversion that acts upon $\mathbb{Z}_2$-graded linear spaces.
A mathematical formulation of the problem can be outlined as follows. We consider a manifold 
\begin{equation}\label{E:fgwqxc}
\mathbf{P}^{n|n+1}=V\times \Pi V\backslash \{0\}/ \mathbb{C}^*
\end{equation} 

 The moduli space of algebraic smooth curves in $\mathbf{P}^{n|n+1}$ of genus $g$ is denoted by  $\M_g(\mathbf{P}^{n|n+1})$. It is a ( nonconnected ) superorbifold.  The problem is to  compute  Berezinian of $\M_g(\mathbf{P}^{n|n+1})$ in simple terms.

The moduli space $\M_g(\mathbf{P}^{n|n+1})$ admits a projection
\begin{equation}\label{E:gsdkskb}
n:\M_g(\mathbf{P}^{n|n+1})\rightarrow {\cal M}_g
\end{equation}
where $\M_g=\M_g(pt)$. The space $\M_g$ carries several line bundles: $\lambda$ is a determinant line bindle. 
It will be described more explicitly  in section (\ref{S:ewecvd}). There is also the canonical line bundle $K=\Omega_{\M_g}^{3g-3}$- a bundle of holomorphic volume forms. There is a classical identity $\lambda^{\otimes 13}=K$  (see \cite{Mumford},\cite{KZ}).
\begin{proposition}\label{P:fgfkbkf}
Denote by $n^*$ a pullback operation by the  map $n$.

If $g=0$, then $\Ber\M_0(\mathbf{P}^{n|n+1})$ is trivial.

For $g=1$ $\Ber\M_1(\mathbf{P}^{n|n+1})$ is equal to $n^*\lambda^{\otimes 2}$.

$\Ber\M_g(\mathbf{P}^{n|n+1})=n^*\lambda^{\otimes 14}$, $g\geq 2$. 
\end{proposition}

The note is organized as follows: in section (\ref{S:ewecvd}) we compute $\Ber\M_g(\mathbf{P}^{n|n+1})$ and prove proposition (\ref{P:fgfkbkf}) using simple homological algebra. In section (\ref{S:jfjhas}) we consider  possible generalizations of the construction from section (\ref{S:ewecvd}).

\section{Computation of $\Ber\M_g(\mathbf{P}^{n|n+1})$}\label{S:ewecvd}

The manifold $\mathbf{P}^{n|n+1}$ splits, the splitting is induced by projection $p:V+\Pi V\rightarrow V$. We use the same notation for projection 
\begin{equation}\label{E:goqdcmmkj}
\mathbf{P}^{n|n+1}\rightarrow \mathbf{P}^{n}
\end{equation}

The following classical  algebro-geometric construction is essential for our computation.

The manifold $\mathbf{P}^{n}$ is a space of complex lines $\{[l]\}$ through the origin in $V$.
 The manifold $\mathbf{P}^{n}$ carries the  Hopf line bundle ${\cal O}(-1)$. Its total space is formed by pairs $\{([l],a)|a\in [l]\subset V\}$. . The total space of ${\cal O}(-1)$ is  embedded as a  subspace into a  direct product $\mathbf{P}^{n}\times V$. The embedding produces a  short exact sequence of bundles
\begin{equation}\label{E:fgsdfg}
0\rightarrow \O \rightarrow V\rightarrow  T(-1)\rightarrow 0
\end{equation}
where $T$ is the  tangent bundle \footnote{It is a standard practice in algebraic geometry to denote a tensor product of sheaves $\F\otimes \O(1)^{\otimes i}, i\in \mathbb{Z}$ by $\F(i)$. An operation of tensor multiplication of $\F$ on $\O(i)$ is called twist.}
of  $\mathbf{P}^{n}$.
It can be easily  verified fiberwise.
  After  appropriate twist the  short exact sequence (\ref{E:fgsdfg}) becomes 
\begin{equation}\label{E:fgs}
0\rightarrow {\cal O} \rightarrow V(1)\rightarrow  T\rightarrow 0
\end{equation}

An alternative description of $\mathbf{P}^{n|n+1}$ is available.
There is an isomorphism $\mathbf{P}^{n|n+1}=\Pi V(1)$. The manifold $\Pi V(1)$ is obtained from the total space of $V(1)$ by the reversal of  parity of the fibers.

This isomorphism can be verified. The quotient (\ref{E:fgwqxc}) can be identified with $V\backslash 0\times \Pi V/ \mathbb{C}^*$ because embedding $V\backslash 0\times \Pi V\subset V\times \Pi V\backslash 0$ is an isomorphism (if we replace $\Pi V$ by some even space , this becomes false ). We can conclude that $V\times \Pi V\backslash \{0\}$ is a total space of a sum of several copies of a line bundle with $\Pi$ applied fiberwise. This line bundle is $V\backslash 0\times \Pi \mathbb{C}/\mathbb{C}^*$, which by definition is equal to $\Pi\O(1)$.

We adopt the following notations : a vector bundle $\F$ over an algebraic manifold $M$ has a space of global sections denoted by  $H^0(M,\F)$ . A linear space of $i$-th cohomology of $\F$ computed through, say, \v{C}ech resolution is denoted by $H^i(M,\F)$. A fibration $p:M\rightarrow N$ enables us to transfer vector bundles from $M$ to $N$ by an operation of direct image $p_*\F$. By definition the fiber of $p_*\F$ at $x\in N$ is equal to $H^0(p^{-1}(x), \F)$. Similarly a fiber of $R^i\F$ at $x$ is $H^i(p^{-1}(x), \F)$ .  A vector bundle $R^i\F$ is called a higher direct image. There is a caveat to the above  definition.  The dimensions of $H^i(p^{-1}(x), \F)$   can change from point to point.  Thus $R^ip_*\F$ is a vector bundle only in favorable circumstances but  typically it is only a sheaf. 

Suppose $M$ is a supermanifold which splits, i.e. there is a fibration  $p:M\rightarrow M_{rd}$ identical on $M_{rd}$  ($M_{rd}\subset M$ is the underlying ordinary manifold). Let $T_M$ be the tangent bundle of $M$. The bundle $T_M$ restricted on $M_{rd}$ splits into a sum of even and odd parts $T^0+T^1$. There is an isomorphism of $T^0$ over $M_{rd}$ with the  tangent bundle  $T_{M_{rd}}$. To simplify notations we drop the subscript $T=T_{M_{rd}}$. If we reverse the  parity of the fibers of $T^1$ we get an ordinary vector bundle $E$ over $M_{rd}$. According to \cite{Manin} the total space of  $\Pi E$ is (noncanonically) isomorphic to $M$ and 
\begin{equation}\label{E:33gdd}
\Ber M\cong p^*(det^{-1}T\otimes det E)
\end{equation}
\begin{definition}
We denote by $det$ the top exterior power of a vector space (bundle).
\end{definition}

Denote by $\M_{g,k}(\mathbf{P}^{n|n+1})$ the moduli space of smooth connected algebraic curves in $\mathbf{P}^{n|n+1}$ with $k$ distinct marked points.
The orbifold $\M_{g,k}(\mathbf{P}^{n|n+1})$ splits over $\M_{g,k}(\mathbf{P}^{n})$:
\begin{equation}
p:\M_{g,k}(\mathbf{P}^{n|n+1}) \rightarrow \M_{g,k}(\mathbf{P}^{n})
\end{equation}
 The splitting comes from the projection (\ref{E:goqdcmmkj}).



 Let us find vector bundles $T$ and $E$ in a  context of the  pair $\M_{g}(\mathbf{P}^{n})\subset \M_{g}(\mathbf{P}^{n|n+1})$. 

A linear space  $T_{\psi}$  is the tangent vector space at a point $\psi\in \M_{g}(\mathbf{P}^{n})$ .
The tangent space to the moduli of curves in a manifold at a point represented by a curve is equal to the space of section of the normal bundle to the curve. 
For $\mathbf{P}^{n}$ it means the following: let $\psi:\Sigma \rightarrow \mathbf{P}^{n}, \psi\in \M_{g}(\mathbf{P}^{n})$.  There is a map $D:T_{\Sigma}\rightarrow \psi^*T_{\mathbf{P}^{n}}$ defined by Jacobian of $\psi$. Denote by $N$ the normal bundle.  It completes the map $D$ to a short exact sequence:
\begin{equation}
0\rightarrow T_{\Sigma}\rightarrow \psi^*T_{\mathbf{P}^{n}} \rightarrow N \rightarrow 0
\end{equation}  

The tangent space  $T_{\psi}$ is equal to $H^0(N)=H^0(\Sigma,N)$\footnote{In the future to simplify notations we suppress $\Sigma$ in $H^i(\Sigma,.)$}. Similarly we define $E_{\psi}$ as $H^0(\psi^*V(1))$. According to the  formula (\ref{E:33gdd}) the fiber of \\ $\Ber\M_{g}(\mathbf{P}^{n|n+1})$ at a point $\psi$ is equal to $p^*(det^{-1}H^0(N)\otimes detH^0(\psi^*V(1)))$.

The  group $H^0(N)$ is a part of a long exact sequence

\begin{equation}
\begin{split}
&0\rightarrow H^0(T_{\Sigma})\rightarrow H^0(\psi^*T_{\mathbf{P}^{n}}) \rightarrow H^0(N) \rightarrow \\
&\rightarrow H^1(T_{\Sigma})\rightarrow H^1(\psi^*T_{\mathbf{P}^{n}}) \rightarrow H^1(N) \rightarrow 0
\end{split}
\end{equation}
It gives an isomorphism of fibers of determinant  line bundles :
\begin{equation}\label{E:jjh22fdd}
\begin{split}
&det^{-1}H^0(N)\otimes det H^0(\psi^*  V(1))=\\
&=det^{-1}H^0(\psi^*T_{\mathbf{P}^{n}})\otimes detH^0(\psi^*  V(1))\otimes det^{-1}H^1(T_{\Sigma})\otimes\\
&\otimes det H^1(\psi^*T_{\mathbf{P}^{n}})\otimes det^{-1}H^1(\psi^*  V(1))\otimes detH^0(T_{\Sigma})
\end{split}
\end{equation}



Machinery of higher direct images is a convenient tool that enables us to see how various vector spaces involved in the formula (\ref{E:jjh22fdd}) vary in families.
A universal curve $\M_{g,1}(\mathbf{P}^{n})$  has evaluation and forgetting maps 
\begin{equation}\label{E:ckdkdk}
\begin{split}
& ev:\M_{g,1}(\mathbf{P}^{n})\rightarrow \mathbf{P}^{n}\\
&q:\M_{g,1}(\mathbf{P}^{n})\rightarrow \M_{g}(\mathbf{P}^{n})
\end{split}
\end{equation}
A fiber of $q$ is equal to a curve $\Sigma$. A vector space $H^0(\psi^*T_{\mathbf{P}^{n}})$ is the fiber at $\psi\in \M_{g}(\mathbf{P}^{n})$ of the bundle $q_*ev^*\mathbf{P}^{n}$. A space $H^0(\psi^*V(1))$ is the fiber of $q_*ev^*V(1)$ at the same point.

A pullback transforms an exact sequence of vector bundles to an exact sequence.
Denote a pullback of the exact sequence (\ref{E:fgs}) by the map $ev$ by
\begin{equation}\label{E:fgwqs}
0\rightarrow \O \rightarrow E_2\rightarrow  E_1\rightarrow 0
\end{equation}
We use the isomorphism $ev^*\O=\O$.

The  direct image functor  by a map  $q$ transforms the short exact sequence (\ref{E:fgwqs})  into a long exact sequence of vector bundles:
\begin{equation}\label{E:jdfsh}
\begin{split}
&0\rightarrow q_*{\cal O} \rightarrow q_*E_2\rightarrow  q_* E_1\rightarrow R^1q_*{\cal O}\rightarrow \\
& \rightarrow R^1q_*E_2\rightarrow  R^1q_* E_1\rightarrow 0
\end{split}
\end{equation}
It generalizes a long exact sequence in cohomology.

Some comments are in order on the sheaves that appear in the last formula.

The fiber of $q_*{\O}$ at $(\Sigma,\psi)$ is  a space of global holomorphic functions over $\Sigma$. It is exhausted by constants for compact $\Sigma$. Thus  $q_*\O=\O$. The fiber of $ R^1q_*{\cal O}$ at the same point is a space dual to the space of holomorphic differentials on $\Sigma$(Serre duality). It does not depend on the map $\psi$; the vector bundle $ R^1q_*{\cal O}$ originates on $\M_g$. This observation will be used later in this note. A fiber of $R^1q_*E_2$ is equal to $H^1(\psi^*V(1))$. The latter is zero  upon using Kodaira vanishing theorem and the fact that   $\O(1)$ is ample. Due to exactness of the sequence (\ref{E:jdfsh}) the sheaf $R^1q_* E_1$  also vanishes. This argument eliminates terms $detH^1(\psi^*T_{\mathbf{P}^{n}})\otimes det^{-1}H^1(\psi^*  V(1))$ in the formula (\ref{E:jjh22fdd}).


Considerations of the previous paragraph imply that the long exact sequence (\ref{E:jdfsh}) gives rise to a canonical isomorphism.
\begin{equation}\label{E:hgsh}
 det(R^1q_*{\cal O})^{-1}= det(q_*E_2)\otimes det^{-1}(q_* E_1)
\end{equation}
The LHS of (\ref{E:hgsh}) is a pullback from $\M_g$.

Let $T_{rel}$ be a line bundle of vector fields tangential to the fibers of projection $q$. The reader can easily identify the fiber of $q_*T_{rel}$ over $\Sigma$ with $H^0(T_\Sigma)$ and the fiber of  $R^1q_*T_{rel}$ with $H^1(T_\Sigma)$.

Suppose $g\geq 2$. Then for any $\Sigma$ $H^0(T_{\Sigma})=0$ and we can rewrite $det^{-1}T\otimes det(E)$ using equations (\ref{E:jjh22fdd}) and (\ref{E:hgsh}) in terms  independent of $\mathbf{P}^{n}$:
\begin{equation}
det^{-1}T\otimes det(E)=det^{-1}(R^1q_*{\cal O})\otimes det^{-1}(R^1q_*T_{rel})
\end{equation}

According to \cite{Mumford} and \cite{KZ} there is an isomorphism between certain line bundles on $\M_{g}$. Let $\lambda=det^{-1}(R^1q_*{\cal O})$ and $K=det^{-1}(R^1q_*T_{rel})$, then $K=\lambda^{\otimes 13}$ and $\Ber\M_{g}(\mathbf{P}^{n|n+1})=\lambda^{\otimes 14}$.

If $g=0,1$, additional terms in the formula for Berezinian (according to isomorphism (\ref{E:jjh22fdd})) are present: $\Ber=det(R^1q_*{\cal O})^{-1}\otimes det^{-1}(R^1q_*T_{rel})\otimes det(R^0q_*T_{rel})$ .

If $g=0$, then $R^1q_*{\cal O}$ and  $R^1q_*T_{rel}$ vanish (direct computation). The group  $H^0(T_{\Sigma})$ (a fiber of $R^0q_*T_{rel}$) is isomorphic to the Lie algebra $\mathfrak{sl}_2$. The group $Aut(\Sigma)$ acts trivially on linear space $det(\mathfrak{sl}_2)$. We conclude that $det(R^1q_*{\cal O})^{-1}\otimes det^{-1}(R^1q_*T_{rel})\otimes det(R^0q_*T_{rel})$ is trivial.

If $g=1$, then  the line bundles  $R^iq_*T^{\otimes k}_{rel}$ are one-dimensional for all $i=0,1$ and $k\in \mathbb{Z}$. Denote $\Omega_{rel}=T^{*}_{rel}$. Using Serre duality we conclude that $$det(R^1q_*{\cal O})^{-1}\otimes det^{-1}(R^1q_*T_{rel})\otimes det(q_*T_{rel})=q_*\Omega_{rel}\otimes q_*\Omega_{rel}^{\otimes 2}\otimes q_*T_{rel}=q_*\Omega_{rel}^{\otimes 2}$$

Thus we have an isomorphism
\begin{equation}
\lambda^{\otimes 2}=det^{-2}R^1q_*\O=q_*\Omega^{\otimes 2}_{ref}=\Ber \M_1(\mathbf{P}^{n|n+1})
\end{equation}
%
%

\section{Computation with targets more general then $\mathbf{P}^{n|n+1}$}\label{S:jfjhas}
In this section we extend proposition (\ref{P:fgfkbkf}) to wider class of manifolds. 

Suppose $M$ is a compact K\"{a}hler manifold. Let $\bigoplus_{i+j=k}H^{i,j}$  be the Hodge decomposition of $k$-th de Rham cohomology; $h^{i,j}=dimH^{i,j}$. Denote $h^{1,1}$ by $h$. We have an identification $H^{1,1}=H^1(M,\Omega^1)$. The last vector space is the Dolbeault cohomology of the vector bundle of  holomorphic differentials. The isomorphisms $H^1(M,\Omega^1)=Ext^1(\O,\Omega^1)=Ext^1(T,\O)$ are standard in algebraic geometry (see \cite{GH}). The  linear space $Ext^1(.,.)$ is a group of extensions of vector bundles. The group $Ext^1(T,\O)$ classifies extensions $0\rightarrow \O \rightarrow X \rightarrow T\rightarrow 0$. Denote by $\mathbb{C}^h$ the group $H^1(M,\Omega^1)$. There is  a  universal extension 
\begin{equation}
0\rightarrow \mathbb{C}^h \rightarrow E \rightarrow T\rightarrow 0
\end{equation}
Denote by $W(M)$ a split supermanifold. It is equal to  the total space of $E$ with the reversed parity of the fibers.
It is convenient to impose some restrictions on $M$:

{\bf Condition} \ \ {\it For any holomorphic curve $\psi:\Sigma \rightarrow M$ a linear space $H^1(\Sigma ,\psi^*E)$ vanishes.}

If $E$ has a filtration with ample adjoint quotients the above condition is fulfilled. One can find many examples of such $M$ in the class of toric varieties. A space $\mathbf{P}^n$ considered in section (\ref{S:ewecvd}) is one of them. We have $h(\mathbf{P}^n)=1$, $E=\mathbb{C}^{n+1}(1)$.
\begin{proposition}
Suppose $M$ satisfies above assumption. The  Berezinian of \\ ${\cal M}_g(W(M))$ is equal to $n^* \lambda^{13+h}$, $g\geq 2$.

If $g=1$, then $\Ber{\cal M}_g(W(M))$ is equal to $n^* \lambda^{h+1}$.

If $g=0$, then $\Ber{\cal M}_g(W(M))$ is trivial.
\end{proposition}
\begin{proof}
The same as for $\mathbf{P}^n$.
\end{proof}
\subsection{Acknowledgements}
This work was written while author was visiting  Institut Mittag-Leffler (2004) and Institute for Advanced Study (2005-2007).
The author wishes to thank these institutions for hospitality. He also would like to thank L. Mason, A.Schwarz and E.Witten  for useful discussions.

\end{document}